\documentclass[10pt,a4paper,fleqn]{article}
\usepackage{amsmath}
\usepackage{amssymb}
\usepackage{amscd}
\usepackage{plain}
\usepackage{graphicx}
\usepackage{epsfig}
\usepackage{fancyhdr}
\usepackage{t1enc,epsfig}
\usepackage[mathscr]{eucal}
\usepackage{epstopdf}
\usepackage{amsthm}

\newtheorem{theorem}{Theorem}
\newtheorem*{thma}{Theorem A}
\newtheorem*{thmb}{Theorem B}
\newtheorem*{thmc}{Theorem C}

\newtheorem{corollary}[theorem]{Corollary}
\newtheorem{definition}[theorem]{Definition}
\newtheorem{example}[theorem]{Example}
\newtheorem{lemma}[theorem]{Lemma}
\newtheorem{proposition}[theorem]{Proposition}
\newtheorem{remark}[theorem]{Remark}
\newtheorem{conjecture}[theorem]{Conjecture}
\begin{document}

\newcommand{\bt}{\begin{theorem}}
\newcommand{\bthma}{\begin{thma}}
\newcommand{\bthmb}{\begin{thmb}}
\newcommand{\bthmc}{\begin{thmc}}
\newcommand{\et}{\end{theorem}}
\newcommand{\ethma}{\end{thma}}
\newcommand{\ethmb}{\end{thmb}}
\newcommand{\ethmc}{\end{thmc}}
\newcommand{\bd}{\begin{definition}}
\newcommand{\ed}{\end{definition}}
\newcommand{\bs}{\begin{proposition}}
\newcommand{\es}{\end{proposition}}
\newcommand{\bp}{\begin{proof}}
\newcommand{\ep}{\end{proof}}
\newcommand{\be}{\begin{equation}}
\newcommand{\ee}{\end{equation}}
\newcommand{\ul}{\underline}
\newcommand{\br}{\begin{remark}}
\newcommand{\er}{\end{remark}}
\newcommand{\bex}{\begin{example}}
\newcommand{\eex}{\end{example}}
\newcommand{\bc}{\begin{corollary}}
\newcommand{\ec}{\end{corollary}}
\newcommand{\bl}{\begin{lemma}}
\newcommand{\el}{\end{lemma}}
\newcommand{\bj}{\begin{conjecture}}
\newcommand{\ej}{\end{conjecture}}
\vskip .5 truecm \def\ux{{\underline x}} \def\umu{{\underline\mu}} \def\omu{{\overline\mu}}

\title{Weight functions and log-optimal investment portfolios}

\author{Y.~Suhov$^{1-3}$, I.~Stuhl$^{4-6}$, M.~Kelbert$^{7,8}$}

\date{}
\maketitle

\begin{abstract}  Following the paper by Algoet--Cover  (1988), we analyse log-optimal portfolios
where return evaluation includes `weights' of different outcomes. The results are twofold: (A) under certain
conditions, logarithmic growth rate is a supermartingale, and (B) the optimal (martingale) investment
strategy is a proportional betting; it does not depend on the form of the weight function, although the
optimal rate does. The existence of an optimal investment strategy has been established earlier in a great
generality by Kramkov--Schachermayer (2003) although our underlying assumptions are different.
\end{abstract}

\footnotetext{$^{1}$ Math Dept, Penn State University, PA, USA;  $^{2}$ DPMMS, University of Cambridge,
 UK;  $^{3}$ IPIT RAS, Moscow, RF}

\footnotetext{$^{4}$ IMS, Univesity of Sao Paulo, SP, Brazil;  $^{5}$ Math Dept, University of Denver, CO, USA; $^{6}$ University of Debrecen, Hungary}

\footnotetext{$^{7}$ Math Dept, University of Swansea, UK;  $^{8}$ Moscow Higher School of Economics,
 RF}

\footnotetext{  \par}

\def\fB{\mathfrak B} \def\fF{\mathfrak F}
\def\fM{\mathfrak M}\def\fX{\mathfrak X}
 \def\cB{\mathcal B}\def\cM{\mathcal M}\def\cX{\mathcal X}
\def\bu{\mathbf u}\def\bv{\mathbf v}\def\bx{\mathbf x}\def\by{\mathbf y}
\def\om{\omega} \def\Om{\Omega}
\def\bbP{\mathbb P} \def\hw{h^{\rm w}} \def\hwi{{h^{\rm w}}}
\def\beq{\begin{eqnarray}} \def\eeq{\end{eqnarray}}
\def\beqq{\begin{eqnarray*}} \def\eeqq{\end{eqnarray*}}
\def\rd{{\rm d}} \def\Dwphi{{D^{\rm w}_\phi}}
\def\BX{\mathbf{X}}\def\Lam{\Lambda}\def\BY{\mathbf{Y}}

\def\mwe{{D^{\rm w}_\phi}}
\def\DwPhi{{D^{\rm w}_\Phi}} \def\iw{i^{\rm w}_{\phi}}
\def\bE{\mathbb{E}}
\def\1{{\mathbf 1}} \def\fB{{\mathfrak B}}  \def\fM{{\mathfrak M}}
\def\diy{\displaystyle} \def\bbE{{\mathbb E}} \def\bu{\mathbf u}
\def\BC{{\mathbf C}} \def\lam{\lambda} \def\bbB{{\mathbb B}}
\def\bbR{{\mathbb R}}\def\bbS{{\mathbb S}} \def\bmu{{\mbox{\boldmath${\mu}$}}}
 \def\bPhi{{\mbox{\boldmath${\Phi}$}}}
\def\bbZ{{\mathbb Z}}
\def\fF{\mathfrak F}
\def\B1{\mathbf 1}
\def\eps{\varepsilon}
\def\ueps{{\underline\eps}}
\def\phi{\varphi}
\def\cF{\mathcal F}

\def\beal{\begin{array}{l}}
\def\beac{\begin{array}{c}}
\def\beacl{\begin{array}{cl}}
\def\ena{\end{array}}
\def\bb{{\mathbf b}}
\def\bpp{{\mathbf p}}
\def\alph{\alpha}
 \def\rA{{\rm A}} \def\vphi{{\varphi}} \def\bear{\begin{array}{r}}
\def\tf{{\tt f}} \def\tg{{\tt g}}

\def\bx{\bar x}

\medskip

{\bf I.} This note is an initial part of a work on log-optimal portfolios influenced by Refs \cite{AC}; see
also \cite{CT}, Chapter 6.
We also intend to use recent progress in studying {\it weighted entropies}; cf. \cite{KM}, \cite{SY} --
\cite{SYS}. A strong impact on the whole direction of research was made by \cite{KrS1}, \cite{KrS2}
where a powerful methodology of a convex analysis have been developed (and elegantly presented)
in a general form, leading
-- among other achievements -- to existence of log-optimal portfolios. See Theorem 1 from \cite{KrS2}.
In the present article, we consider a situation of an arbitrary weight (or utility) function which does not fall
under assumptions imposed in \cite{KrS2}. Moreover, we go beyond existence and provide a specific form
of the optimal strategy.

The result offered here is as follows.

You are betting on results $\eps_n$ of subsequent random trials, $n=1,2,\ldots$. Each $\eps_n$ produces
a value $x_n\in\cX_n$ where $(\cX_n,\fX_n,\mu_n)$ is assumed to be a standard measure
space. We suppose that a random string $\ueps_1^n=\left(\beac\eps_1\\ \vdots \\ \eps_n\ena\right)$
has a joint probability
density function (PDF) $f_n(\ux_1^n)$ relative to reference measures $\omu_n=\prod\limits_{j=1}^n\mu_j$,
on $\operatornamewithlimits{\times}\limits_{j=1}^n\fX_j$:
\beq\label{eq1}
\bbP (\ueps_1^n\in A) =\int_Af_n(\ux_1^n)\rd\omu_n (\ux_1^n),\;\ux_1^n=\left(\beac x_1\\ \vdots\\ x_n\ena
\right)\in\operatornamewithlimits{\times}\limits_{j=1}^n\cX_j,\,
A\in\operatornamewithlimits{\times}\limits_{j=1}^n\fX_j .
\eeq
A conditional PDF $\tf_n(x_n|\ux_1^{n-1})$ will be also used, with
\beq\label{eq2}
\tf_n(x_n|\ux_1^{n-1})f_{n-1}(\ux_1^{n-1})=f_n(\ux_1^n)\hbox{ and }\;
\int_{\cX_n}\tf_n(x_n|\ux_1^{n-1})\rd\mu_n(x_n)=1,\,\hbox{a.s.}
\eeq

Let us agree that if you stake $\$\,C_n$ on game $n$ you win $\$\,C_ng_n(x_n)$ if the result is
$x_n\in\cX_n$. (So, you make a profit when $C_ng_n(x_n)>0$ and incur a loss when $C_n
g_n(x_n)<0$.) Here $g_n$ are given real-valued functions $x_n\in\cX_n\mapsto g_n(x_n)\in
\bbR$.$^{*)}$\footnote{$^{*)}
$All functions figuring throughout the paper are assumed measurable, with a specific indication of the
sigma-algebra when necessary.}  We say that $g_n$ are return functions.
%It is supposed that $\int_{\cX_n}{\mathbf 1}(g(x_n)<0)\rd\mu_n(x_n) >0$.
\def\ueps{{\underline\eps}}

Let $Z_0>0$ be an initial capital. More generally, given $n\geq 1$, denote by $Z_{n-1}>0$ your
fortune after $n-1$ trials and impose the restriction that variable
$C_n=C_n(\ueps_1^{n-1})$ is $\fF_{n-1}$-measurable. Here and below, $\fF_0=\sigma (Z_0)$
and $\fF_n=\fF_0\vee
\left(\operatornamewithlimits{\times}\limits_{j=1}^n\fX_j\right)$ for $n\geq 1$. (One says that $C_n$
is a previsible strategy.)  Then $Z_{n-1}=Z_{n-1}
(\ueps_1^{n-1})$ is $\fF_{n-1}$-measurable. It also makes sense to require that
$C_n\geq 0$.$^{**)}$\footnote{$^{**)}$ One
 also may demand that $-C_ng_n(x_n)\leq Z_{n-1}$ for
$\mu_n$-a.a. $x_n\in\cX_n$. (In applications, this is required to guarantee the deposit.)}
 We have the recursion
\beq\label{eq3}
Z_n=Z_{n-1}+C_ng_n(\eps_n)=Z_{n-1}\left(1+\frac{C_ng_n(\eps_n)}{Z_{n-1}}\right)
\eeq
and wish to maximize $\bbE S_N$ where
\beq\label{eq4}
S_N:=\sum\limits_{j=1}^N\phi_j(\eps_j;\ueps_1^{j-1})\log\,\frac{Z_j}{Z_{j-1}}.
\eeq
Here the weight function (WF) $\ux_1^j\mapsto \phi_j(x_j;\ux_1^{j-1})\geq 0$ depends on $\ux_j$
and the vector $\ux_1^{j-1}$.  Quantity
$\phi_j(x_j;\ux_1^{j-1})$ represents a `sentimental' value of outcome $x_n$ (given that it succeeds
a sequence $\ux_1^{j-1}$) taken into account when one calculates $S_N$.  Value $\bbE S_N$
is the weighted expected interest rate after $N$ rounds of investment.
When $\phi_j\equiv 1$, the sum \eqref{eq4} becomes telescopic and equal to $\diy\log\,\frac{Z_N}{Z_0}$,
the standard interest rate. Recursion \eqref{eq3} suggests a martingale-based approach.

We also consider a sequence of positive functions $b_n(x_n)$,
%=\diy\frac{1}{(2\pi)^{d/2}}\exp\,\left(-\frac{\|\bx\|^2}{2}\right)$,
$x_n\in\cX_n$, figuring  in Eqns \eqref{eq5} -- \eqref{eq7}. More precisely, we will use
the following conditions \eqref{eq5}, \eqref{eq6}.
\beq\label{eq5}
\int_{\cX_n}\phi_n (x_n ;\ueps_1^{n-1})b_n(x_n)g_n(x_n)\rd\mu_n (x_n) =0,\;\hbox{a.s.}
\eeq

\def\rA{{\rm A}}

\beq\label{eq6}
\beal\diy\int_{\cX_n}\phi_n(x_n ;\ueps_1^{n-1})b_n (x_n )\rd\mu_n(x_n)
\leq \bbE\big[\vphi_n(\eps_n;\ueps_1^{n-1})\big|\fF_{n-1}\big]\\
\qquad{}\qquad{}\qquad{}\diy =\int_{\cX_n}
\phi_n(x_n;\ueps_1^{n-1})\tf_n(x_n|\ueps_1^{n-1})\rd \mu_n(x_n),\;\hbox{a.s.}\ena
\eeq
Next, define a RV $\alph_n=\alph_n(\ueps_1^{n-1})$ by
\beq\label{eq7}
\beacl\diy\alph_n&\diy =\bbE\;\,\left\{\vphi_n(\eps_n;\ueps_1^{n-1})\log\,\frac{\tf_n(\eps_n|
\ueps_1^{n-1})}{b_n (\eps_n)}\Big|\fF_{n-1}\right\}\\
\;&\diy=\int_{\cX_n}\phi_n(x_n;\ueps_1^{n-1})\tf_n(x_n|\ueps_1^{n-1})
\log\,\frac{\tf_n(x_n|\ueps_1^{n-1})}{b_n(x_n)}
\rd\mu_n(x_n).\ena
\eeq

\bthma \label{thm:A} Given $1< N\leq\infty$,
assume that functions $\vphi_n$ are non-negative and obey
conditions \eqref{eq5}, \eqref{eq6} for $1\leq n<N$. Then:

{\rm{(a)}} For all previsible $C_n$ such that $\diy1+\frac{C_ng(\eps_n)}{Z_{n-1}}>0$,
sequence $S_n-\rA_n\;$ is a supermartingale, where
$\rA_n:=\diy\sum\limits_{j=1}^n\alpha_j$. Consequently,
$\bbE\,S_n\leq \diy\sum\limits_{j=1}^n\bbE\,\alpha_j$.

{\rm{(b)}}  Sequence $S_n-\rA_n$, $1\leq n<N$,
is a martingale for some previsible $C_n$ satisfying $0\leq C_n\leq Z_{n-1}$ and
$\diy1+\frac{C_ng_n(\eps_n)}{Z_{n-1}}>0$ a.s.
iff the following holds. There exists a function $\ux_1^{n-1}\mapsto D_{n-1} (\ux_1^{n-1})\in [0,1]$ with
$1+D_n(\ueps_1^{n-1})g_n(\eps_n)>0$ a.s. such that
\beq\label{eq8}
\tf_n(x_n|\ux_1^{n-1})=g_n(x_n)b_n(x_n)D_{n-1} (\ux_1^{n-1})+b_n(x_n),\,\hbox{a.s.}
\eeq
In this case
%a non-negative value not depending on $x_n\in\cX_n$, and
\beq\label{eq9}
C_n (\ueps_1^{n-1})=D_{n-1}(\ueps_1^{n-1})Z_{n-1}(\ueps_1^{n-1}).
\eeq
\ethma

\bp (a) Write:
$$\beal\bbE\Big\{\big(S_n-\rA_n\big)\big|\fF_{n-1}\Big\}=S_{n-1}-\rA_{n-1}\\
\qquad\quad\diy+\bbE\left\{\left[\phi_n (\eps_n;\ueps_1^{n-1})\log\,\left(1+\frac{C_ng_n(\eps_n)}{Z_{n-1}}
\right)\right]\Big|
\fF_{n-1}\right\}-\alph_n.\ena$$
Next, represent
\beq\label{eq10}
\beal\diy\bbE\left\{\left[\phi_n (\eps_n;\ueps_1^{n-1})\log\,\left(1+\frac{C_ng_n(\eps_n)}{Z_{n-1}}\right)
\right]\Big|
\fF_{n-1}\right\}-\alph_n \\
\quad\diy =\;\int_{\cX_n}\; \phi_n(x_n;\ueps_1^{n-1})\;\tf_n(x_n|\ueps_1^{n-1})\;\log\;\left[1 +
\frac{C_ng_n(x_n)}{Z_{n-1}}\right] \rd\mu_n(x_n)\\
\qquad\qquad\diy - \int_{\cX_n}\phi_n(x_n;\ueps_1^{n-1})\tf_n(x_n|\ueps_1^{n-1})\log\,
\frac{\tf_n(x_n|\ueps_1^{n-1})}{b_n (x_n)}\rd\mu_n(x_n)\\
\quad\diy = \int_{\cX_n}\phi_n(x_n ;\ueps_1^{n-1})\tf_n(x_n|\ueps_1^{n-1}) \log\frac{1
+C_ng_n(x_n)/Z_{n-1}}{\tf_n(x_n|\ueps_1^{n-1})/b_n (x_n )} \rd\mu_n(x_n)\\
\quad\diy = \int_{\cX_n}\phi_n(x_n ;\ueps_1^{n-1})\;\tf_n(x_n|\ueps_1^{n-1})\;
\log\;\frac{h_n(x_n)}{\tf_n(x_n|\ueps_1^{n-1} )}\; \rd\mu_n(x_n)\\
\quad\diy \leq \int_{\cX_n}\phi_n(x_n;\ueps_1^{n-1})\tf_n(x_n|\ueps_1^{n-1})\\
\qquad\qquad\qquad\times\diy\;\left[ \frac{h_n(x_n )}{\tf_n
(x_n|\ueps_1^{n-1})} - 1 \right]{\mathbf 1}\big(\tf_n
(x_n|\ueps_1^{n-1})>0\big)\; \rd\mu_n(x_n)\\
\quad\diy = \int_{\cX_n}\phi_n(x_n ;\ueps_1^{n-1}) \Big[h_n(x_n) - \tf_n(x_n|\ueps_1^{n-1})
\Big]\rd\mu_n(x_n) \leq 0,
\ena
\eeq
where $h_n(x_n) := \diy b_n (x_n )\left[1 + \frac{C_ng_n(x_n )}{Z_{n-1}}\right]$, $x_n\in\bbR^d$.
The final inequality in \eqref{eq10} holds since, almost surely,
\beq\label{eq11}
\beal
\diy \int_{\cX_n}\phi_n(x_n ;\ueps_1^{n-1})h_n(x_n ) \rd\mu_n(x_n) = \int_{\cX_n}
\phi_n(x_n ;\ueps_1^{n-1})b_n (x_n ) \rd\mu_n(x_n)\\
\qquad\qquad\diy + \int_{\cX_n}\phi_n(x_n ;\ueps_1^{n-1})b_n (x_n)\frac{C_ng_n(x_n)}{Z_{n-1}}
\rd\mu_n(x_n)\\
\qquad\qquad\qquad\qquad\qquad\qquad\diy \leq \int_{\cX_n}\phi_n(x_n;\ueps_1^{n-1})
\tf_n(x_n|\ueps_1^{n-1})
\rd\mu_n(x_n),\ena
\eeq
due to \eqref{eq5} and \eqref{eq6}.

As a result, we get the supermartingale inequality
\beq\label{eq12}
\bbE\left\{\Big[S_n-\rA_n\Big]\Big|\fF_{n-1}\right\}\leq S_{n-1}-
\rA_{n-1},\;\hbox{ a.s.}
\eeq

(b) For the martingale property we need to fulfill equalities in Eqn \eqref{eq10}. The first inequality
becomes equality iff $\diy\left[ \frac{h_n(x_n )}{\tf_n(x_n|\ueps_1^{n-1})} - 1 \right]{\mathbf 1}\big(\tf_n
(x_n|\ueps_1^{n-1})>0\big)=0$ $\mu_n$-a.s., i.e.,
$$b_n (x_n )\left[1 + \frac{C_ng_n(x_n )}{Z_{n-1}}\right]=\tf_n(x_n|\ueps_1^{n-1}),\;\;
\tf_n(\,\cdot\,|\ueps_1^{n-1})\;\hbox{-a.s.,}$$
which is Eqn \eqref{eq8}. The second inequality in \eqref{eq11} also follows from \eqref{eq8}.
\ep

\medskip
{\bf Remarks.} {\bf 1.} The martingale strategy, when it exists, provides a log-optimal
investment portfolio.

{\bf 2.} Quantities $g_n(x_n)$ and $b_n(x_n)$ can be made dependent on
argument $\ux_1^{n-1}$ as well; in this case $h_n(x_n)$ also becomes a function of $x_n$
and $\ux_1^{n-1}$. In fact, functions $(x_n;\ux_1^{n-1})\mapsto g_n(x_n;\ueps_1^{n-1})$ can be
considered as a part of the investment strategy. Taking $b_n(x_n;\ux_1^{n-1})=
\tf_n(x_n|\ux_1^{n-1})$ leads to the (non-interesting) case $S_n=0$.

{\bf  3.} The staple of the proof of Theorem A \ref{thm:A} is the Gibbs inequality for weighted
entropies; see  \cite{KM}, \cite{SY} -- \cite{SYS}. It is similar to the standard Gibbs inequality (cf. \cite{CT},
\cite{KS}) but requires additional assumptions, as listed in Theorem A \ref{thm:A}.

{\bf  4.} In fact, the inequality in Eqn \eqref{eq10} may hold when $\vphi_n$ is not necessarily
non-negative; in such a situation, methods of  convex analysis developed and used in \cite{KrS1}, \cite{KrS2}
would not be suitable. However, even assuming that functions $\vphi_n \geq 0$, our conditions in
Theorem A\ref{thm:A} cover a variety of cases left open by Theorem 1 from \cite{KrS2}.
At the same time, the Gibbs inequality can be considered as a special fact from convex analysis;
thus, connections between our methodology and the one from \cite{KrS2} need further explorations.

Another feature of Theorem A\ref{thm:A} is that it specifies an optimal policy.

\medskip

{\bf II.} The level of generality adopted in Theorem A\ref{thm:A} may seem excessive from the point
of view of applications. We therefore provide a special form of the statement where
trials $\eps_n$ are IID, and each trial produces one of $m>1$ outcomes
$E_1$,$\ldots$, $E_m\in\bbR$ with probabilities  $p_1$,$\ldots$, $p_m >0$. We also set
the return function $g_n(E_i)=E_i$ and use
uniform probabilities to emulate functions $b_n$: $b_n(E_i)=\diy\frac{1}{m}$. Here if you stake
$\$\,C_n$ on game $n$ you win $\$\,C_nE_i$ if the result is $E_i$.
As above, let $Z_{n-1}>0$ the fortune
after $n-1$ trials ($Z_0>0$ is the initial capital).  As before, let $\fF_n =\sigma (Z_0)$
and $\fF_n =\sigma (Z_0,\ueps_1^n)$, $n\geq 1$, and consider a sequence  of RVs $C_n$
where $C_n$ is $\fF_{n-1}$-measurable (a previsible strategy).
The recursion \eqref{eq3} becomes
\beq\label{eq13}
Z_n=Z_{n-1}+C_n\eps_n=Z_{n-1}\left(1+\frac{\eps_nC_n}{Z_{n-1}}\right).
\eeq

\noindent
We wish to maximize, in $C_n$, the weighted expected interest rate $\bbE S_N$ where
\beq\label{eq14}
S_n:=\sum\limits_{j=1}^n\phi (\eps_j)\log\,\frac{Z_j}{Z_{j-1}}.
\eeq
Here $E\mapsto\phi (E)\geq 0$ is a
weight function (for simplicity depending only upon a one-time outcome).

Theorem A\ref{thm:A} then takes the following form:

\bthmb\label{thm:B} Suppose that
\beq\label{eq15}
\sum\limits_i\phi (E_i)E_i=0\;\hbox{ and }\;\;
\frac{1}{m}\sum\limits_i\phi (E_i)\leq \sum\limits_i\phi (E_i)p_i.
\eeq
Set:
\beq\label{eq17}
\alpha =
\sum\limits_i\phi (E_i)p_i\log\,(p_im).
\eeq
Then

{\rm{(a)}} For all previsible $C_n$ with $\diy 1+\frac{\eps_nC_n}{Z_{n-1}}>0$,
sequence $S_n-\alpha n\;$ is a supermartingale; consequently,
$\bbE\,S_n\leq n\alph$.

{\rm{(b)}} $S_n-\alpha n$ is a martingale for a previsible $C_n$ with
$0\leq C_n\leq Z_{n-1}$ and $\diy 1+\frac{\eps_nC_n}{Z_{n-1}}>0$ iff $\diy D:=\frac{mp_i-1}{E_i}$
is a non-negative number between $0$ and $1$ which does not depend on outcome $E_i$, and
$C_n=DZ_n$.
\ethmb

In case $m=2$, the above martingale strategy exists only if $E_1=-E_2$
and $\phi (E_1)=\phi (E_2)$ (no weight preference). Assume for definiteness that $E_1>0$
and $p_1\geq 1/2$. Then $D=\diy\frac{2p_1-1}{E_1}=\frac{1-2p_1}{E_2}$, and the martingale strategy is
$\diy C_n=\frac{Z_{n-1}}{E_1}(2p_1-1)$.
It means that you repeatedly bet the proportion
$\diy\frac{2p_1-1}{E_1}$ of your current capital on outcome $E_1$.
\medskip

{\bf III.} Another example of interest is where $\cX_n=\bbR^d$ and $\mu_n$ is a standard Lebesgue's
measure. Setting
\beq\label{eq19}
f_n(x_n|\ux_1^{n-1})=\frac{\exp\Big(-x_n^{\rm T}\Sigma^{-1}x_n/2\Big)}{[(2\pi )^d{\rm det}\,\Sigma ]^{1/2}},
\;x_n\in\bbR^d,\;\ux_1^{n-1}\in (\bbR^d)^{n-1},
\eeq
yields IID Gaussian random vectors $\eps_n\sim{\rm N}(0,\Sigma)$. Let us take
\beq\label{eq20}
b_n(x_n)=\frac{\exp\Big(-x_n^{\rm T}\Sigma_0^{-1}x_n/2\Big)}{[(2\pi )^d{\rm det}\,\Sigma _0]^{1/2}},
\eeq
where $\Sigma_0\neq\Sigma$. Also let us fix a return function $x\in\bbR^d\mapsto g(x)$ and consider  a
weight function
$x\in\bbR^d\mapsto \varphi (x)$ depending on the current outcome $x$ only. Then Theorem A
\ref{thm:A} transforms into Theorem C\ref{thm:C}:

\bthmc\label{thm:C}
Suppose that $\phi$ and $g$ satisfy
\beq\label{eq21}
\int_{\bbR^d}\phi (x)g(x)\exp\Big\{-\frac{1}{2}\Big[x^{\rm T}\Big(\Sigma^{-1}+\Sigma_0^{-1}\Big)x
\Big]\Big\}\rd x =0
\eeq
and
\beq\label{eq22}
\frac{1}{[{\rm det}\,\Sigma_0]^{1/2}}\int_{\bbR^d}e^{-x^{\rm T}\Sigma_0^{-1}x/2}\phi (x)\rd x\leq
\frac{1}{[{\rm det}\,\Sigma]^{1/2}}\int_{\bbR^d}
\phi ( x )e^{-x^{\rm T}\Sigma^{-1}x/2}\rd x.\quad\eeq
Define value $\alph$ by
\beq\label{eq23}\beal\diy
\alph =\frac{1}{[(2\pi)^d{\rm det}\Sigma]^{1/2}}\int_{\bbR^d}\phi ( x)e^{- x^{\rm T}\Sigma^{-1} x/2}\\
\qquad\qquad\diy\times \frac{1}{2}\Big\{x^{\rm T}\Big(\Sigma_0^{-1} -\Sigma^{-1}\Big)x+\log\Big[{\rm det}\,\big(\Sigma_0\Sigma^{-1}\big)\Big]\Big\}\rd x.\end{array}\eeq
Then:

{\rm{(a)}} For all previsible $C_n$ with $\diy1+\frac{C_ng (\eps_n)}{Z_{n-1}}>0$, sequence
 $S_n-n\alph\;$ is a supermartingale, and hence
$\bbE\,S_n\leq n\alph$.

{\rm{(b)}} $S_n-n\alph$
is a martingale for some previsible $C_n$ with $0\leq C_n\leq Z_{n-1}$ and
$\diy1+\frac{C_ng (\eps_n)}{Z_{n-1}}>0$ a.s. iff, for some constant $D\in (0,1)$ the
strategy is $C_n (\ueps_1^{n-1})=DZ_{n-1}(\ueps_1^{n-1})$ and return function
$g(x)$, $x\in\bbR^d$, has the form
%a non-negative value not depending on $x_n\in\cX_n$, and
\beq\label{eq24}
g(x)=\frac{1}{D}\left\{[{\rm det}(\Sigma_0\Sigma^{-1})]^{1/2}e^{- x^T\Big(\Sigma^{-1}-\Sigma_0^{-1}\Big)
x/2} -1\right\}.\eeq
\ethmc

{\bf Remark.} The statement of Theorem C\ref{thm:C} can be repeated for any choice of
two PDFs in Eqns \eqref{eq19}, \eqref{eq20}, with an obvious modification of \eqref{eq24}.

\vskip 10 pt

\noindent
{\bf Acknowledgement}
\vskip 10pt
The authors are grateful to their respective institutions for the financial support and hospitality during
the academic year 2014-5: YS thanks the Math Dept, PSU; IS was supported by the
FAPESP Grant No 11/51845-5, and thanks the IMS, USP and to Math Dept, University of Denver; MK thanks the Higher School of Economics for the support in the framework of the Global Competitiveness Programme.

\end{document}